\documentclass[runningheads]{llncs}
\usepackage{graphicx}
\usepackage{amsmath,amssymb,bbm}
\usepackage{latexsym,todonotes}
\usepackage{mathtools}
\usepackage{cite}
\usepackage{color}
\usepackage{arydshln}
\usepackage{subfigure}
\usepackage{bm}
\usepackage{enumitem}
\usepackage{mathrsfs}
\usepackage{epstopdf}
\usepackage{algorithmicx,algorithm}
\usepackage{algpseudocode}
\usepackage[bookmarks=false]{hyperref}
\usepackage{tikz,xcolor}

\usepackage{textcomp}
\def\BibTeX{{\rm B\kern-.05em{\sc i\kern-.025em b}\kern-.08em
    T\kern-.1667em\lower.7ex\hbox{E}\kern-.125emX}}

\graphicspath{{figs/}}

\definecolor{lime}{HTML}{A6CE39}
\DeclareRobustCommand{\orcidicon}{%
	\begin{tikzpicture}
	\draw[lime, fill=lime] (0,0) 
	circle [radius=0.16] 
	node[white] {{\fontfamily{qag}\selectfont \tiny ID}};
	\draw[white, fill=white] (-0.0625,0.095) 
	circle [radius=0.007];
	\end{tikzpicture}
	\hspace{-2mm}
}

\foreach \x in {A, ..., Z}{%
	\expandafter\xdef\csname orcid\x\endcsname{\noexpand\href{https://orcid.org/\csname orcidauthor\x\endcsname}{\noexpand\orcidicon}}
}

\mathcode`*=\string"8000
\begingroup
\catcode`*=\active
\xdef*{\noexpand\textup{\string*}}
\endgroup
\begin{document}
\title{Energy Efficient Priority-Based Task Scheduling for Computation Offloading in Fog Computing\thanks{The work described in this paper is supported by Student Interdisciplinary Research Fund from BNU-HKBU United International College, College Research Grant from BNU-HKBU United International College R201911, and Zhuhai Basic and Applied Basic Research Foundation Grant ZH22017003200018PWC.}}
\titlerunning{Priority-Based Policy for Offloaded Task Scheduling in Fog Computing}
%
\author{Jiaying Yin\inst{1}\orcidA{} \and
Jing Fu\inst{2}\orcidB{}\and
Jingjin Wu\inst{1}\orcidC{} \and Shiming Zheng\inst{1}\orcidD{}}
\authorrunning{J. Yin et al.}
%
\institute{Division of Science and Technology, BNU-HKBU United International College, Zhuhai, Guangdong, P. R. China \\
\and
School of Engineering, RMIT University, Victoria, Australia \\
\email{n830006230@mail.uic.edu.cn, jing.fu@rmit.edu.au, jj.wu@ieee.org, l630005069@alumni.uic.edu.cn}}
\maketitle              
\setcounter{footnote}{0}
\begin{abstract}
Fog computing offers a flexible solution for computational offloading for Internet of Things (IoT) services at the edge of wireless networks. It serves as a complement to traditional cloud computing, which is not cost-efficient for most offloaded tasks in IoT applications involving small-to-medium levels of computing tasks. Given the heterogeneity of tasks and resources in fog computing, it is vital to offload each task to an appropriate destination to fully utilize the potential benefit of this promising technology. In this paper, we propose a scalable priority-based index policy, referred to as the Prioritized Incremental Energy Rate (PIER), to optimize the energy efficiency of the network. We demonstrate that PIER is asymptotically optimal in a special case applicable for local areas with high volumes of homogeneous offloaded tasks and exponentially distributed task durations. In more general cases with statistically different offloaded tasks, we further demonstrate the improvement of PIER over benchmark policies in terms of energy efficiency and the robustness of PIER to different task duration distributions by extensive simulations. Our results show that PIER can perform better than benchmark policies in more than $78.6\%$ of all simulation runs.

\end{abstract}

\section{Introduction}\label{sec:Introduction}

Fog computing, a term originally introduced by Cisco, is considered a refinement to the traditional cloud computing where intensive computational tasks are offloaded to the data centers in the central cloud by migrating part of cloud computing and storage capabilities to the edge of the network~\cite{ciscofog}. This concept is particularly suitable for the realization of the Internet of Things (IoT). In IoT applications such as smart cities, remote gaming, and instant healthcare, offloading intensive computational tasks from mobile terminals (MTs) is frequently required because of the limited computing power of the MTs. In the IoT era, {\it fog nodes} (FNs) at the network edge, including base stations, micro data centers, or other access devices, are usually considered more appropriate offloading destinations as the central cloud is often geographically far from the MTs. The high latency caused by propagation between the MT and the central cloud would be unacceptable for these delay-sensitive applications~\cite{Mukherjee2018}.

Another advantage of fog computing is higher energy efficiency~\cite{Jalali2016}. The MTs are usually operated by batteries with a limited lifespan and thus are not suitable for most computational tasks. At the same time, the data centers in the central cloud are generally much less energy efficient for offloading tasks of small to moderate sizes. Furthermore, additional power consumption is required for long-distance transmissions between the MT and the central cloud. The FNs, closer to the MTs and connected to the electrical grid, offer a reliable platform for energy-efficient computations.

While FNs have higher computing and storage capabilities than MTs, the types of resources and the number of resources for each kind in an FN are still limited compared to the central cloud~\cite{Jalali2016}. Meanwhile, the latency and the energy consumption for a particular task generated by an MT may vary if it is offloaded to different FNs. Two main factors causing the difference are the quality of the wireless connection between the MT and the FN and the availability of required resources in the FN. Therefore, from the perspective of an {\it Edge Infrastructure Provider} (EIP) who owns and operates the FNs, an effective task scheduling policy is essential to fully exploit the benefits of fog computing.

In this paper, we consider a heterogeneous fog computing architecture with different classes of tasks. The task classes are differentiated by the resources required for each task. At the same time, FNs have distinct energy profiles and different sets of resources that can execute the offloaded computing tasks from MTs. In addition, the allocation of offloading tasks is not only subject to the availability of computing resources in FNs but also constrained by the number of wireless channels between transmitting MT and FN. On the other hand, the central cloud can be considered as the last resort of offloading when all computing resources in FNs are occupied, as additional latency and power consumption would be incurred if a task is offloaded further. Our objective in this paper is to propose a task scheduling policy aiming to optimize the energy efficiency, namely minimizing the ratio of long-run average power consumption of the network to the long-run average throughput, subject to the availability of transmission channels and computing resources.

Existing work on computational offloading policies in fog computing mainly adopts a static optimization approach by considering a single instance of the network or assuming that the arrival pattern of tasks is fully known. It has been shown that, under such conditions, an efficient offloading policy could improve utilization of resources in FNs~\cite{Deng2016}, reduce the energy consumption of the network~\cite{Jalali2016}, guarantee fairness among MTs in terms of transmission rate~\cite{Du2018}, improve the Quality of Experience of mobile users~\cite{He2018}, or achieve an appropriate trade-off between the QoS requirement of users and energy-related cost of the EIP~\cite{You2016}. While static optimization approaches can restrict the complexity of the problem to a reasonable level, the main shortcoming of such methods is that the stochastic nature of network traffic is ignored. As most offloaded tasks by personal users have relatively short durations, the availability of computational resources is expected to change frequently over time. As a result, static optimization approaches are likely to miss the potential benefit gained by dynamically reusing the resources released by completed tasks. While more recent studies (e.g.~\cite{Chang2020, Gao2021}) proposed online algorithms that consider the stochastic factor, various important issues, such as the heterogeneity of MTs and FNs, concurrent constraints from a limited number of channels and resources, and non-linear cost functions, are still ignored. Such deviations may affect the practicability of proposed policies.

Motivated by the facts above, we propose a scalable and robust \emph{priority-based index policy} in this paper, which makes real-time decisions on the offloading destination based on the instantaneous state of the network upon the initiation of an offloading task to maximize the energy efficiency of the network. We will demonstrate that the performance gap between the proposed policy and the optimal solution diminishes as the system size increases in a particular case appropriate for a small-area environment where the offloaded tasks are generated from densely distributed mobile devices, while both tasks and FNs are largely homogeneous (e.g., large sports events~\cite{ZHANG2020}). Such property is desirable for the problem studied in this paper, as offloading in fog computing generally involves a large-scale system with huge numbers of MTs and FNs. In more general cases, we will also demonstrate that the proposed policy outperforms benchmark policies adopted in existing work. 

While similar ideas have been applied in our previous work~\cite{fu2018restless, Wang2018, Fu2020}, we now consider a tailor-made model for the offloading problem in fog computing, where a task may be offloaded to an FN at the network edge or the central cloud, subject to the availability of both computing resources and wireless transmission channels. We will also consider that the power consumption at the network edge is, in fact, a piece-wise function of the load, as an "idle power" would be needed to activate the hardware equipment supporting a resource group when at least one unit of the resources is occupied~\cite{Jalali2016}. To summarize, the model we consider in this paper involves many practical features in task scheduling in fog computing and thus is significantly different from the fundamental models in existing studies.

We will describe the model in Section \ref{sec:model} and the stochastic optimization problem in Section \ref{sec:problem}.
We then introduce the proposed scheduling policy in Section \ref{sec:policy}, present numerical results in Section \ref{sec:simulation}, and give conclusions in Section \ref{sec:conclusion}.

\section{Network Model}\label{sec:model}

Let $\mathbb{R}_0$, $\mathbb{R}_+$ and $\mathbb{N}_+$ represent the sets of non-negative reals, positive reals and positive integers, respectively. Let $[N]$ represent a set $\{1,2,\ldots,N\}$ for any $N\in\mathbb{N}_+$.

We consider a wireless network with orthogonal transmission channels~\footnote{Orthogonal wireless channel allocations eliminate intra-cell interference and utilize frequency spectrum resources more efficiently.}, through which wireless connections between mobile terminals (MTs) and fog nodes (FNs) at the network edge are established. Computing tasks generated by MTs can be appropriately offloaded to FNs through the transmission channels. For the sake of presentation, we refer to all these storage and computing resources in FNs and the cloud as \emph{abstracted resource components} (ARCs), which are classified in \emph{groups} (referred as \emph{ARC groups} thereafter) based on their functionalities and geographical locations. Denote $K$ as the number of ARC groups at the network edge, and let $C_k\in\mathbb{N}_+$ represent the \emph{capacity} of ARC groups $k\in[K]$, which is the total number of ARC units supporting the execution of offloaded tasks in the ARC group. 

FNs are connected to the central cloud through wired cables, where tasks can be further offloaded. As the capacity of the cloud is much larger than the capacity of any ARC groups $k\in[K]$ in FNs~\cite{Du2018}, we assume that there are always sufficient resources in the cloud to support all types of offloaded tasks, that is, the cloud has infinite capacity. 

Similarly, consider $J\in\mathbb{N}_+$ classes of offloaded tasks, classified by the differences in their locations of origin, application styles, and relevant requirements on computing resources.
As mentioned earlier, an offloaded task can be completed at the network edge or further offloaded to the cloud. We refer to the tasks belonging to class $j\in[J]$ as $j$-tasks. Tasks of the same type can be those generated by MTs in the same area, with the same application styles, and has the same requirements on ARCs. 

If a $j$-task is accommodated by ARC group $k\in[K]$, $w_{j,k}\in[C_k]$ units of ARCs will be occupied by this task.
These occupied units will be released upon the completion of the task (that is, when the computational results are transmitted back to the respective MT) and can be reused by future tasks~\cite{Gillam2018}.
While, if a $j$-task cannot be served by ARC group $k$, because of mismatch in geographical locations or required functions, we set $w_{j,k} \rightarrow +\infty$ to prohibit $j$-tasks from being assigned there. 

We divide the network edge into $L\in\mathbb{N}_+$ \emph{destination areas}, and denote $\mathscr{K}_{\ell}\subset [K]$ as the set of ARC groups in destination area $\ell \in [L]$. We assume, without loss of generality, that all $\mathscr{K}_{\ell}$ are non-empty and mutually exclusive for different $\ell\in[L]$.
For each sender-receiver pair implied by MTs generating $j$-tasks and ARC groups in destination area $\ell$, the transmission rate is \emph{i.i.d.} with a mean of $\mu_{j,\ell} \in \mathbb{R}_+$. We further denote $\Gamma_{j,\ell} \in\mathbb{N}_+$ as the number of wireless transmission channels for a $(j,\ell)$ pair (referred as $(j,\ell)$-channels hereafter), and assume that each $j$-task offloaded to destination area $\ell$ occupies one $(j,\ell)$-channel. As transmission to the cloud must be via the network edge, a $j$-task offloaded to the cloud may choose an arbitrary destination area $\ell \in [L]$ and occupy one corresponding $(j,\ell)$-channel. By this definition, it is impossible for a $j$-task to be offloaded when all $(j,\ell)$-channels for all $\ell \in [L]$ are occupied by existing tasks. 

We consider the situation where the computing powers at the network edge or in the cloud are sufficient to complete the tasks in a relatively short period. In this sense, the durations of computational operations of offloaded tasks are generally much shorter than the transmission time; that is, the processing time is dominated by the transmission time between the MT and the FN for offloaded tasks completed at the network edge, and transmission time between the MT and the central cloud via the network edge for tasks offloaded to the central cloud~\cite{You2016, Du2018}. 

By the arguments above, the expected duration of $j$-tasks computed by an ARC group $k\in\mathscr{K}_{\ell}$ at the network edge are assumed to take independently and identically distributed random values, with a mean of $1/\mu_{j,\ell}$.

For offloading to the cloud, two transmission segments, namely from the MT to the network edge and from the network edge to the cloud, are needed. For notation consistency, we define a particular ARC group $k_c$, which has an infinite capacity $C_{k_c} \to +\infty$. As a cable backbone network usually connects the FNs and the cloud, we can assume that the transmission time between the network edge and the cloud, denoted by $D_0$, is the same for all tasks via all destination areas. As the number of channels is limited, tasks offloaded to the cloud will always transmit via the area equipped with the fastest available channel. 
Therefore, the effective mean service rate for tasks offloaded to the cloud via destination area $\ell$ is
\begin{equation*}\label{eqn:cloud_duration}
 \mu'_{j, \ell} = \frac{1}{\left(\frac{1}{\mu_{j,\ell}} + D_0\right)}.
\end{equation*}

We assume that the arrivals of tasks generated from user class $j\in[J]$ form a Poisson process with a mean rate $\lambda_j\in\mathbb{R}_+$. Power consumption of the network edge consists of \emph{idle} and \emph{operational} power. The idle power is the basic consumption incurred when an ARC group is activated (which means the hardware components supporting the ARC group must be powered on). At the same time, the operational power is consumed when ARC units are engaged in processing tasks~\cite{8475536}. The relevant components of ARCs, deployed in FNs such as micro BSs or embedded servers, can be de-activated when not in use for power saving, and activated upon the arrival of a task occupying at least one ARC unit~\cite{Wu2020}\footnote{While the more general cases with non-negligible activation delay and power consumption can be addressed by integrating the vacation queuing model with activation cost and delay as in~\cite{Wu2020}, they will complicate the analysis and we do not consider them in this paper due to the limited space.}. In a dynamic system such as a wireless network, the amount of the operational power consumption is affected by the real-time loads of ARC groups.
We let $\varepsilon_k^0\in\mathbb{R}_+$ and $\varepsilon_k\in\mathbb{R}_+$ represent the amounts of the idle power and operational power consumption per computing unit of ARC group $k\in[K]$. Similar power consumption models have been empirically justified and widely applied in existing research, where idle power and operational power values depend on specific hardware~\cite{Jalali2016, Wu2020}.
Detailed discussion about the power consumption of the network edge will be provided in Section~\ref{sec:problem}.
Similarly, we denote the power consumption of an ARC unit for computing a $j$-task in the cloud as $\bar{\varepsilon}_j\in\mathbb{R}_+$. As demonstrated in existing research, a larger amount of power consumption is required for transmission when a task is offloaded further~\cite{Du2018}. Therefore, for all $j\in[J]$ and $k\in[K]$, we have $w_{j,k}\varepsilon_k <\bar{\varepsilon}_j$ by definition.

Finally, we introduce a \emph{scaling parameter} $h\in \mathbb{N}_+$ to emphasize the scale of the network. We define $\lambda_j \coloneqq h\lambda_j^0$
for all $j\in[J]$, $\varepsilon^0_k \coloneqq h\tilde{\varepsilon}^0_k$, $C_k \coloneqq hC^0_k$ for all $k\in[K]$, and $\Gamma_{j, \ell} \coloneqq h\Gamma_{j, \ell}^0$ for all $j\in[J], \ell\in[L]$ where $\lambda_j^0, \tilde{\varepsilon}^0_k, C^0_k, \Gamma_{j, \ell}^0\in\mathbb{R}_+$. The assumption that $\lambda_j$, $C_k$, $\varepsilon^0_k$ and $\Gamma_{j, \ell}^0$ all increase proportionally with $h$ is reasonable in practical scenarios in edge computing, that the capacity of wireless channels and ARC groups are designed just to meet the potential demand. Correspondingly, a large $h$ is more appropriate for modelling the network environment in a densely populated area where tasks are generated more frequently.

\section{Stochastic Optimization Problem}\label{sec:problem}

Define $[K]^* = [K]\cup\{k_c\}$ as the set of ARC groups in the network edge and the cloud.
Denote $X_{j,k}(t)\in\mathbb{N}_0$, where $j\in[J]$, $k\in[K]^*$, as the number of $j$-tasks being served by ARC units in group $k$ at time $t\geq 0$. Also, we let $Y_{j, \ell}(t)$ represents the number of occupied $(j,\ell)$-channels at time $t$.
Because of the limited capacities of ARC groups and finite number of transmission channels, these variables should satisfy
\begin{equation}\label{eqn:capacity_constraint:1}
\sum\limits_{j\in[J]} w_{j, k}X_{j, k}(t) \leq C_k,~\forall k\in[K],~t \geq 0,
\end{equation}
\begin{equation}\label{eqn:channel_constraint}
\sum\limits_{k\in\mathscr{K}_{\ell}}X_{j, k}(t) \leq \Gamma_{j,\ell},~\forall j \in [J],~\forall \ell\in[L],~t \geq 0,
\end{equation}
\begin{equation}\label{eqn:total_channel_constraint}
\sum\limits_{k\in[K]^*}X_{j, k}(t) \leq \sum\limits_{\ell \in [L]}\Gamma_{j,\ell},~\forall j \in [J], ~t \geq 0.
\end{equation}

Let $\boldsymbol{X}(t)=(X_{j, k}(t): j\in[J], k\in[K]^*$), and the state space of process $\{\boldsymbol{X}(t),\ t\geq 0\}$ (the set involves all possible values of $\boldsymbol{X}(t)$ for $t\geq 0$) be
\begin{equation}
\mathscr{X}=\prod\limits_{j\in[J]}\prod\limits_{k\in[K]^*}\left\{0,1,\ldots,\mathcal{A}_k \right\},
\end{equation}
where $\prod$ represents Cartesian product and $\mathcal{A}_k$ represents an upper limit of the maximum number of requests that can be concurrently accommodated by ARC group $k \in [K]^*$. Specifically,

\begin{equation*}
\mathcal{A}_k=
\begin{cases}
\min\left\{\lfloor \frac{C_k}{w_{j,k}}\rfloor, \Gamma_{j, \ell(k)}\right\}, & \text{if } k \in [K],\\
\sum\limits_{j \in [J]}\sum\limits_{\ell \in [L]}\Gamma_{j,\ell}, & \text{otherwise.}
\end{cases}
\end{equation*}

Note that although $\mathscr{X}$ is larger than the set of possible values of $\boldsymbol{X}(t)$, the process $\boldsymbol{X}(t)$ will be further constrained by~\eqref{eqn:capacity_constraint:1} to~\eqref{eqn:total_channel_constraint} in our optimization problem that will be defined later in this section.

Similarly, we define action variables $a_{j, k}(\boldsymbol{x})\in\{0,1\}$, $\boldsymbol{x}\in\mathscr{X}$, as a function of the state space for each $j\in[J]$ and $k\in[K]^*$.
If $a_{j, k}(\boldsymbol{x})=1$, then $w_{j,k}$ units of ARC in group $k$ are selected to serve an incoming $j$-task when $\boldsymbol{X}(t)=\boldsymbol{x}$; and otherwise, ARC group $k$ is not selected.
Specifically, an incoming $j$-task when $\boldsymbol{X}(t)=\boldsymbol{x}$ means the first $j$-task coming after and exclude time $t$; that is, the process $\boldsymbol{X}(t)$ is defined as left continuous in $t\geq 0$.

We define the action space, which is the set involving all possible values of action variables, as
\begin{equation}
\mathscr{A} = \{0,1\}^{J\times (K+1)}.
\end{equation}

Let $\boldsymbol{a}(\boldsymbol{x})=(a_{j, k}(\boldsymbol{x}): j\in[J],k\in[K]^*)$, $\boldsymbol{x}\in\mathscr{X}$.
For any $\boldsymbol{a}(\boldsymbol{X}(t)) \in \mathscr{A}$ and $\boldsymbol{X}(t)\in\mathscr{X}$, it should further satisfy
\begin{equation}\label{eqn:action_constraint}
\sum\limits_{k\in[K]^*}a_{j, k}(\boldsymbol{X}(t)) \leq 1,~\forall j\in[J].
\end{equation}
Since the probability of tasks simultaneously generated by users in the same class is assumed to be zero, it is reasonable to always select a single ARC group to serve an incoming $j$-task as described by \eqref{eqn:action_constraint}.

A scheduling policy is determined by the action variables for all states $\boldsymbol{x}\in\mathscr{X}$.
We add a superscript and rewrite the action variables as $a^{\phi}_{j, k}(\boldsymbol{x})$ representing the action variables under policy $\phi$,
which can then be considered as a mapping from $\mathscr{X}$ to $\mathscr{A}$.
Similarly, since the stochastic process $\{\boldsymbol{X}(t),t\geq 0\}$ is conditioned on the underlying policy, we rewrite previously defined state variable $\boldsymbol{X}(t)$ as $\boldsymbol{X}^{\phi}(t)$. We let $\Phi$ represents the set of all such policies $\phi$.

The long-run average throughput and power consumption of the network system under a policy $\phi$ are given by \eqref{eqn:throughput} and \eqref{eqn:consumption}, respectively. 

\begin{figure*}
\begin{equation}\label{eqn:throughput}
\begin{split}
\mathcal{L}^\phi \coloneqq \lim\limits_{T \to \infty} \frac{1}{T} \left[\sum\limits_{j\in[J]}\sum\limits_{k\in[K]} \int_0^T \mu_{j,\ell(k)} X^{\phi}_{j, k}(t)\ dt \right. \\
\left. +\sum\limits_{j\in[J]}\sum\limits_{\ell\in[L]}\int_0^T \mu'_{j, \ell} \left( Y_{j, \ell}(t)-\sum\limits_{k \in \mathscr{K}_{\ell}}X^{\phi}_{j, k}(t) \right)\ dt \right]
\end{split}
\end{equation}
\end{figure*}

\begin{figure*}

\begin{equation}\label{eqn:consumption}
\begin{split}
\mathcal{E}^\phi \coloneqq \lim\limits_{T \to \infty} \left[\sum\limits_{j\in[J]}\sum\limits_{k\in[K]}\left( \frac{1}{T} \int_0^T \varepsilon_k w_{j, k}X^{\phi}_{j, k}(t)\ dt\right)
+\sum\limits_{k\in[K]}\left(\frac{1}{T}\int_0^T \varepsilon_k^0 \Theta(X^{\phi}_{j, k}(t))\ dt \right) \right. \\
\left. + \frac{1}{T} \sum\limits_{j\in[J]}\int_0^T \left[\sum\limits_{\ell\in[L]}\left(Y_{j, \ell}(t)- \sum\limits_{k\in\mathscr{K}_{\ell}} X^{\phi}_{j, k}(t)\right)\ dt\right]\overline{\varepsilon}_j \right]
\end{split}
\end{equation}
\end{figure*}

The two terms in \eqref{eqn:throughput} represent the long-run average throughput of the network edge and the central cloud, respectively. For \eqref{eqn:consumption}, the first and second terms represent the long-run average operational and idle power at the network edge, the last term stands for the long-run average power consumption for computing tasks offloaded to the cloud. In addition, $\ell(k)$ represents the destination area $\ell\in[L]$ for $k\in\mathscr{K}_{\ell}$, and $\Theta(x)$ for $x\in\mathbb{R}$ is the Heaviside function, defined as
\begin{equation*}
\Theta(x)=
\begin{cases}
1, & \text{if } x >0,\\
0, & \text{otherwise.}
\end{cases}
\end{equation*}

\section{Scheduling Policy}\label{sec:policy}

We aim to optimize the \emph{energy efficiency} of the network, which is equivalent to minimizing the ratio of the long-run average power consumption to the long-run average throughput:

\begin{equation}\label{eqn:obj}
\min\limits_{\phi \in\Phi}\mathcal{E}^\phi/{\mathcal{L}^\phi},
\end{equation}

We refer to the optimization problem with the objective function~\eqref{eqn:obj} and constraints~\eqref{eqn:capacity_constraint:1} to~\eqref{eqn:total_channel_constraint} and \eqref{eqn:action_constraint} as the \emph{task offloading scheduling problem} (TOSP). The TOSP forms an instance of the \emph{resource allocation problem} discussed in \cite{fu2018restless}, which consists of parallel \emph{restless multi-armed bandit problems} (RMABPs) coupled by capacity constraints. A policy $\phi$ is considered better if it produces a smaller value of the objective function~\eqref{eqn:obj}.

\subsection{Index Policy}

Conventional dynamic optimization techniques are generally not applicable for the TOSP as the computational complexity is prohibitively high for fog computing systems of practical scales due to the large state space of $\boldsymbol{X}(t)$. Therefore, we consider an \emph{index policy}, which only requires a given sequence of \emph{indices} to be assigned to all possible processes $\boldsymbol{X}^{\phi}(t)$ for a certain policy $\phi \in\Phi$. The indices will then determine the priorities of ARC groups to allocate relevant resources for an incoming task.

We propose an offloading policy, called \emph{Prioritized Incremental Energy Rate} (PIER), and intuitively explain the computation of indices in PIER as follows.

When an incoming $j$-task is offloaded to the network edge at time $t$, the EIP needs to decide which ARC group to allocate for handling the task. Note that the power consumption of the network will be increased no matter where the task is computed. We refer such increment as the \emph{incremental power rate} $c_{j, k}(t)$ hereafter. Specifically, suppose the task is computed at the network edge with resources from ARC group $k\in [K]$. In that case, the incremental power rate is $\varepsilon^0_k + w_{j,k}\varepsilon_k$ if none of the ARC units in the group $k$ are occupied upon the arrival of the task, or $w_{j,k}\varepsilon_k$ otherwise. The expected incremental power rate for computing a $j$-task in the cloud is $\bar{\varepsilon}_j$. We compute the indices of ARC group $k$ for a $j$-task, as the ratio of the service rate of ARC group $k$ (determined by the destination area that the group is in) to the incremental power rate of the ARC group if selected, given by $\mu_{j, \ell(k)}/c_{j, k}(t)$. Heuristically, the PIER ranks all ARC groups (including those at the network edge and central cloud) with sufficient resource units by descending order of the indices for a task and then always selects the best ARC group at the moment to optimize the energy efficiency.

By comparing the index values of all feasible offloading destinations for task $j \in [J]$ at time $t$, when the network state is $\boldsymbol{X}(t) \in \mathscr{X}$, we can obtain the set of available ARC groups with the highest index, denoted as,
\begin{equation}\label{eqn:index_state}
\tilde{\mathcal{V}}_j(\boldsymbol{X}(t))\coloneqq \arg\max_{k \in \mathcal{K}_j(t)} [\mu_{j,\ell(k)} / c_{j, k}(t)],
\end{equation}
where $\mathcal{K}_j(t)$
is the set of ARC groups that satisfy the channel and capacity constraints~\eqref{eqn:capacity_constraint:1} to~\eqref{eqn:total_channel_constraint} for an incoming $j$-task at $t$. 
Note that the task will be served if $|\tilde{\mathcal{V}}_j(\boldsymbol{X}(t))|$ is no less than one, and will be blocked if $|\tilde{\mathcal{V}}_j(\boldsymbol{X}(t))| = 0$.
In cases where $|\tilde{\mathcal{V}}_j(\boldsymbol{X}(t))| > 1$,  more than one ARC groups are as efficient as each other, and ties can be broken arbitrarily. We then define that, for any $\boldsymbol{a}(\boldsymbol{X}(t)) \in \mathscr{A}$ and $\boldsymbol{X}(t)\in\mathscr{X}$,
\begin{equation}
a_{j, k}^{\text{PIER}}(\boldsymbol{X}(t))=\left\{\begin{matrix}
1, & \quad \text{if $k = \min \left[\tilde{\mathcal{V}}_j(\boldsymbol{X}(t))\right]$,} \\
0, & \quad \text{otherwise,}
\end{matrix}\right.
\end{equation}
where $\min[\tilde{\mathcal{V}}_j(\boldsymbol{\cdot})]$ returns the minimal element in set $\tilde{\mathcal{V}}_j(\boldsymbol{\cdot})$ given that $\tilde{\mathcal{V}}_j(\boldsymbol{x}) \ne \varnothing$. When $\tilde{\mathcal{V}}_j(\boldsymbol{x})=\varnothing$, $a_{j, k}^{\text{PIER}}(\boldsymbol{X}(t))$ is set to $0$ for all $k$, and the incoming $j$-task is considered blocked. Furthermore, for an incoming $j$-task, the computational complexity of PIER is linear in the number of ARC groups that the task has access to, which is capped at $K+1$. In addition, only binary information that whether each ARC group has extra capacity to accommodate the incoming task is required for making a decision.

\subsection{Asymptotic Optimality}

As in~\cite[Corollary 1]{fu2018restless}, an index policy that prioritizes the process with the highest indices is \emph{asymptotically optimal} as $h \to \infty$. The proposed PIER, while specifically designed for the TOSP, is in the same vein as priority-based policies proposed in~\cite{fu2018restless, Wang2018, Fu2020}, which have been shown to approach the optimal solution in stochastic optimization problems as the scale of the system is sufficiently large. Particularly, when $J = 1$, $K = L$ (that is, there is only one ARC group in every destination area at the network edge), $\Gamma_{j, \ell}^0 = 1$ for all $j\in[J],$ $\ell\in[L]$, $C^0_k = 1$ for all $k \in [K]$, and the duration of tasks are exponentially distributed, the PIER is equivalent to a special case of a policy that has been proposed and proved to be asymptotic optimal in~\cite{Fu2020}. We will demonstrate the asymptotic optimality of the PIER with a numerical example in Section~\ref{sec:simulation}. For more general cases, due to the space constraint, we will also numerically demonstrate and compare the performances of PIER and other benchmark policies in the TOSP in Section~\ref{sec:simulation}.

\section{Numerical Results}\label{sec:simulation}

In addition to the proposed PIER policy, we consider two benchmark policies for performance comparison: the \emph{Prioritized Transmission Rate} (PTR) policy and the \emph{Prioritized Least Power Consumption} (PLPC) policy. Under the PTR policy, an incoming task always offloads to the ARC group with the fastest transmission rate, given that the capacity constraint is satisfied, to minimize the user's waiting time for task completion. The PLPC policy, on the other hand, always chooses the ARC group with the minimal additional energy consumption rate among all groups with available capacity for an incoming task to reduce the energy bill for EIPs.

For all the simulation results presented in this paper, the 95\% confidence intervals based on the Student t-distribution are within 5\% of the observed mean.
\subsection{Verification of Asymptotic Optimality}
We first consider a particular example with a single task class ($J = 1$), five destination areas with a single accessible ARC group in each destination area ($K = L = 5$), each $(j, \ell)$ pair has the equivalent number of wireless transmission channels ($\Gamma_{j, \ell}^0 = 1$), identical capacities across all ARC groups ($C_k^0 = 1$ for all $k \in [K]$), negligible idle power ($\varepsilon^0_k = 0$) and relatively high arrival rates such that some tasks will have to be offloaded to the cloud or dropped as the computing and communication resources at the network edge are insufficient to accommodate all requests. This scenario is appropriate for fog computing in densely populated environments, such as significant sports events, where tasks and available resources are largely homogeneous, and the relevant hardware components are always kept active (thus, the idle power consumption under all policies is the same and can be ignored for comparison purposes) to satisfy the enormous offloading demands~\cite{ZHANG2020}. Under these conditions, the PIER is a special instance of the policy that has been proved to be asymptotically optimal in~\cite{Fu2020}. That is, as the scaling parameter $h \to \infty$, the energy efficiency of PIER approaches the optimal solution.

\begin{figure}
\centering
\includegraphics[height=5.0cm,width=7.5cm]{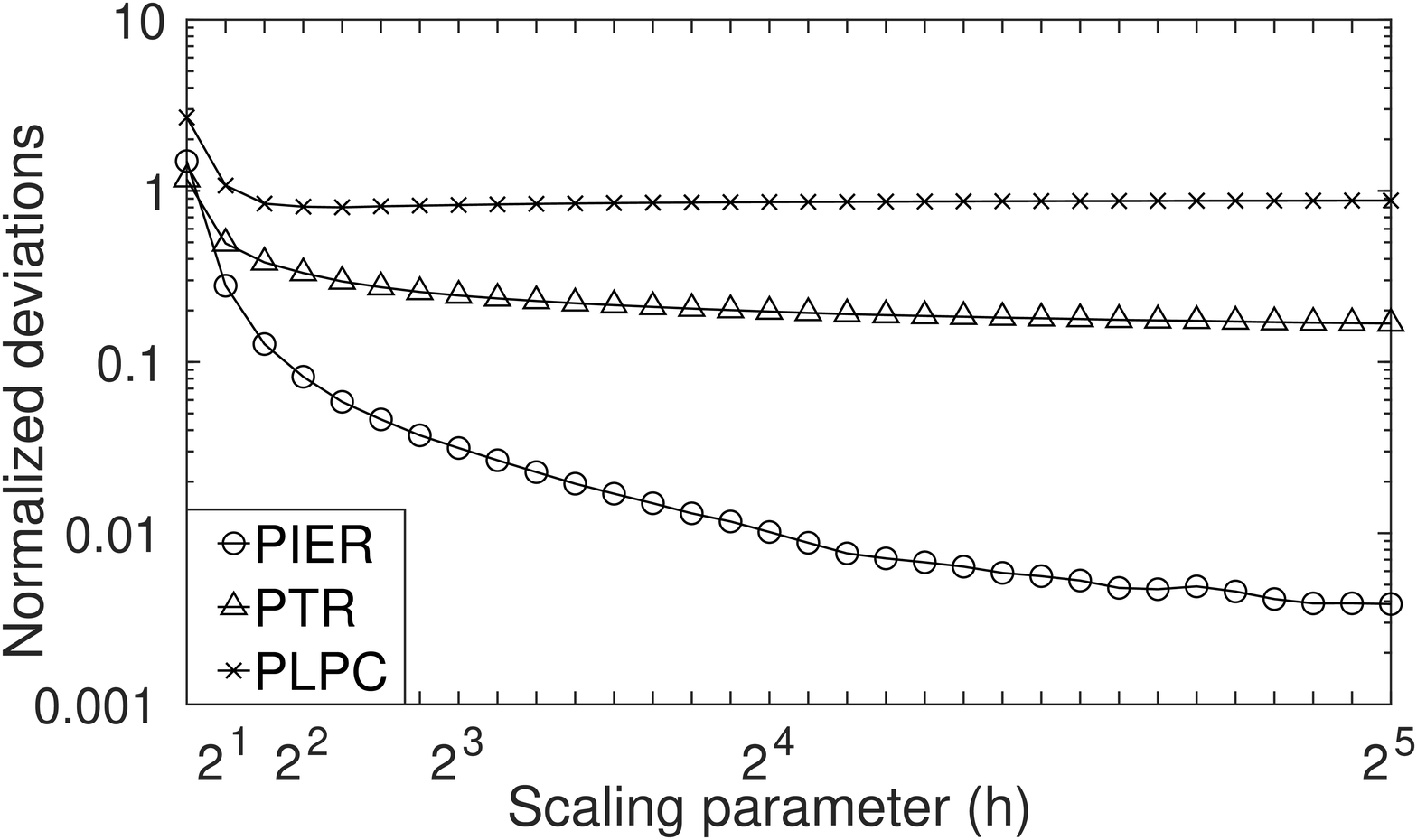}

\caption{Normalized deviations of PIER and benchmark policies versus scaling parameter $h$.\label{fig:fig1}}
\end{figure}

The above conclusion is verified numerically in Fig.~\ref{fig:fig1}, where we compare the energy efficiency of PIER and two benchmark policies for different values of $h$. The durations of tasks served in each destination area are exponentially distributed with mean values $0.587051$, $2.65982$, $0.547387$, $1.1986949$, and $4.78274$, respectively, and the fixed transmission time from the edge to the cloud is $D_0 = 5$. The power consumption per unit of occupied computing resources $\varepsilon_k$ are set to be $1.08316$, $10.0584$, $1.18651$, $8.0544$, and $16.0324$ for each of the five ARC groups at the network edge, and the power consumption for completing one task at the cloud is $\bar{\varepsilon}_j = 51.4714$.
Also, we set the $\lambda_j^0$ for the only task class to be $5.182638$ and consider $w_{j,k} = 1$ for all $k \in [K]$. The \emph{normalized deviation} in Fig.~\ref{fig:fig1} is defined as
$$ \frac{\mathcal{E}^{\phi}/\mathcal{L}^{\phi}-\mathcal{E}^{\phi^O}/\mathcal{L}^{\phi^O}}{\mathcal{E}^{\phi^O}/\mathcal{L}^{\phi^O}}$$
for policy $\phi$, where $\phi^O$ represents the optimal policy. As shown in Fig.~\ref{fig:fig1}, the normalized deviation for PIER is below $1\%$ for $h \ge 10$ and converges to 0 as $h$ increases further, while the normalized deviations of PTR and PLPC stay relatively high even for large values of $h$. Therefore, in this case, PIER is the most appropriate policy, due to its significant superior performance in terms of energy efficiency compared to PTR and PLPC, and advantage in terms of computational complexity over the optimal solution.

\subsection{Performance Comparison in General Cases}

We now consider more general cases, with two task classes and five ARC groups at the network edge, where the number of computing resources units required for each class in all ARC groups ($w_{1,k}$ or $w_{2,k}$ for all $k \in K$) is randomly generated from $\{1, 2\}$. The five ARC groups are distributed in four destination areas, all of which are accessible for both task classes. Tasks from $j$-class could be transmitted to destination area $\ell$ through a ($j, \ell$)-channel, where the number of ($j, \ell$)-channels for each ($j, \ell$) pair is randomly generated from $\{6, 7, 8\}$. For each ARC group, the capacity $C_k$ is randomly generated from $\{3,4,5,6\}$, and the power consumption per unit of occupied resources is randomly generated from $[0.1,20]$. We also set the idle power consumption $\varepsilon^0_k = 0.5\varepsilon_k C_k$ according to the power consumption profile of a small scale FN such as a micro cellular gateway~\cite{Jalali2016}, and keep $D_0 = 5$.

Note that the computational complexity to obtain the optimal solution, as demonstrated in~\cite{Wang2018,Fu2020}, is extremely high even in the special case demonstrated in the previous subsection. It is computationally prohibitive to obtain the optimal solution in the general cases considered in this subsection. Therefore, we will focus on the improvement of PIER over the benchmark policies. 

\begin{figure*}[ht]
\centering
\subfigure[]{\includegraphics[height = 3cm, width=0.32\linewidth]{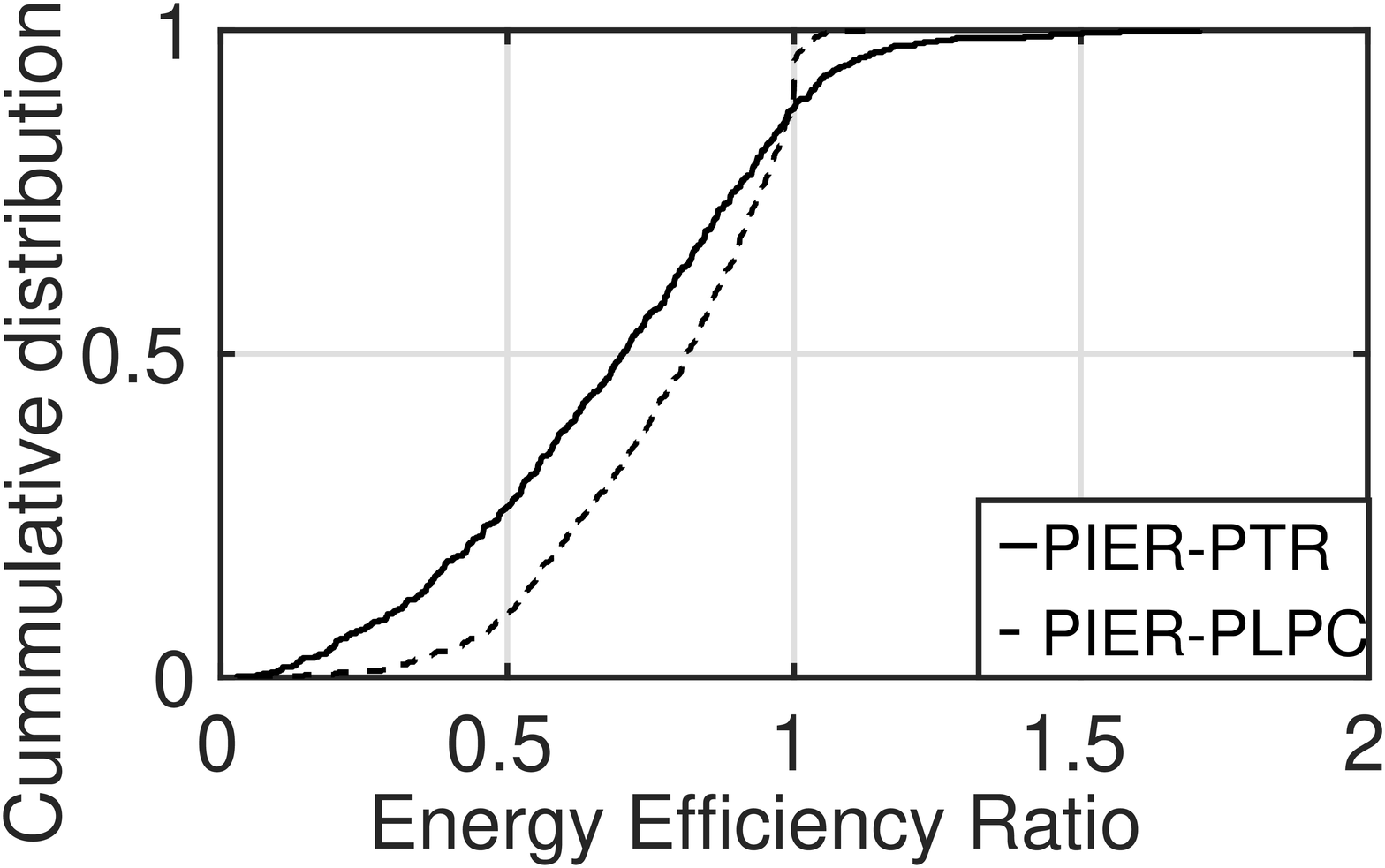}\label{fig:fig2a}}
\subfigure[]{\includegraphics[height = 2.95cm, width=0.32\linewidth]{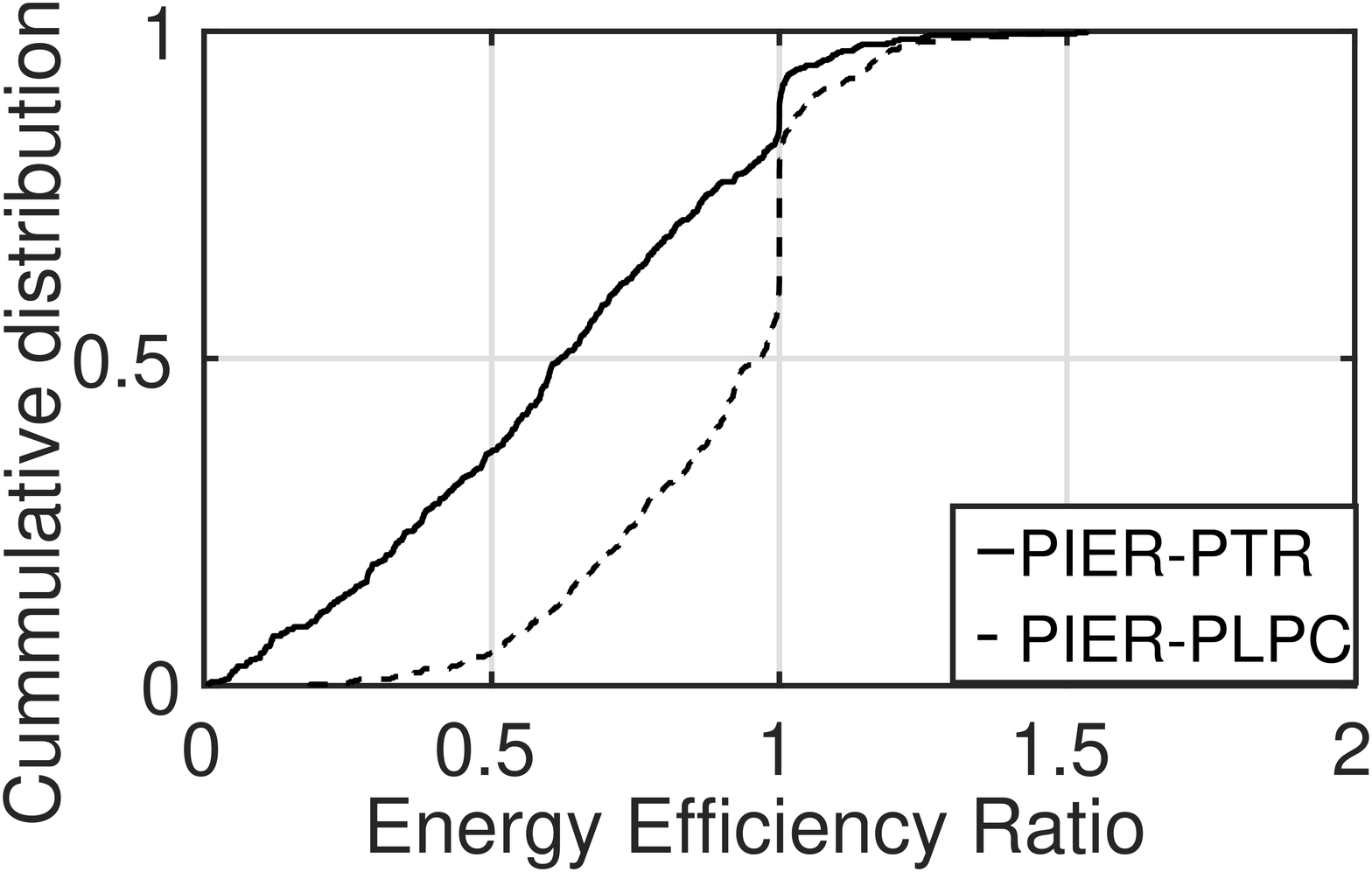}\label{fig:fig2b}}
\subfigure[]{\includegraphics[height = 2.9cm, width=0.32\linewidth]{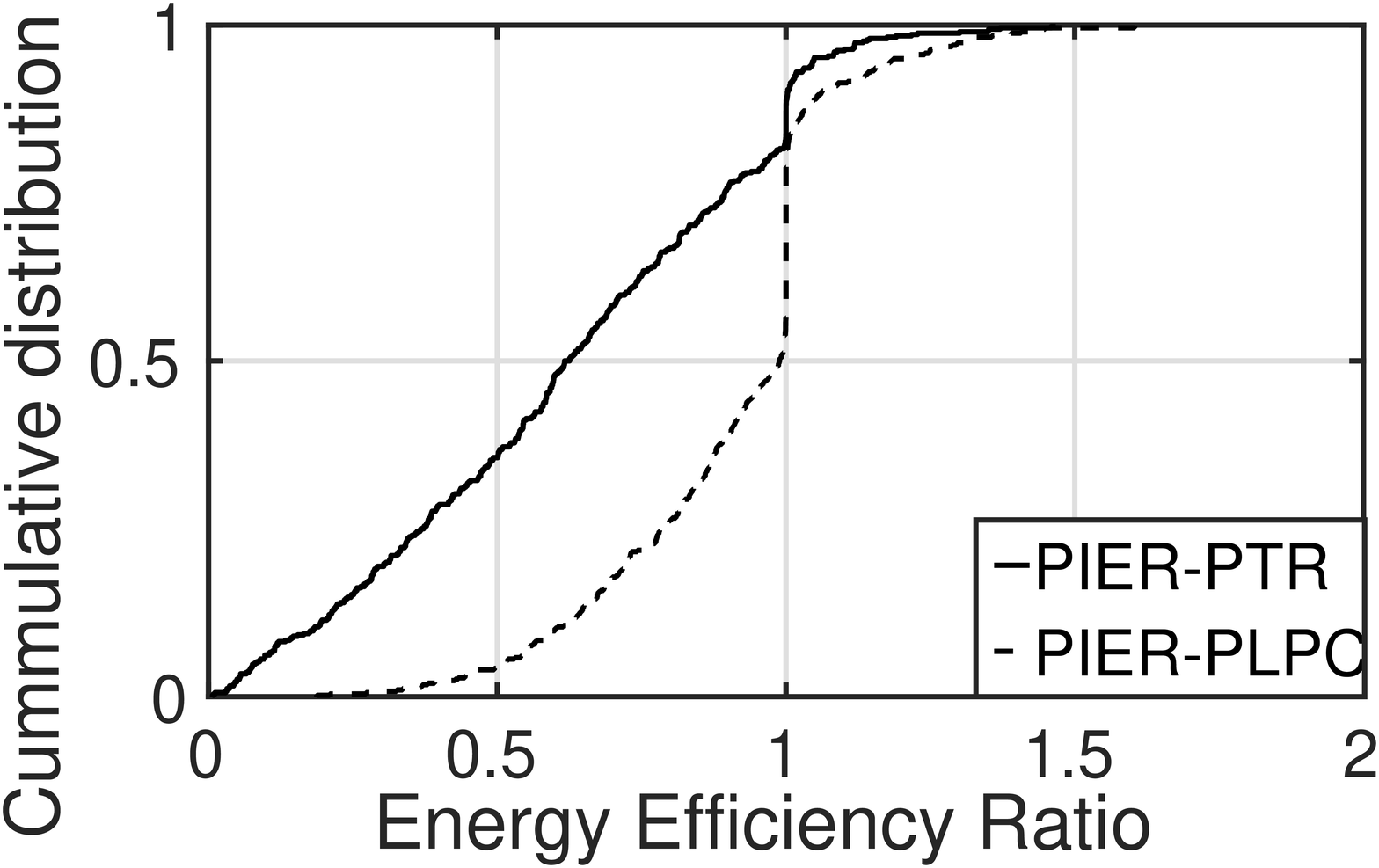}\label{fig:fig2c}}
\caption{Cumulative distribution of energy efficiency ratios of PIER to that of the benchmark policies PTR and PLPC with different scaling parameters: (a) $h=1$; (b) $h=10$; (c) $h=20$.}\label{fig:fig2}
\end{figure*}

In Fig.~\ref{fig:fig2}, we demonstrate the improvement of energy efficiency achieved by PIER as compared to the two benchmark policies in aforementioned scenarios. The cumulative distribution curves of energy efficiency ratios ($\displaystyle \frac{\mathcal{E}^{\phi_1}/\mathcal{L}^{\phi_1}}{\mathcal{E}^{\phi_2}/\mathcal{L}^{\phi_2}}$, where $\phi_1$ and $\phi_2$ represent two different policies, a ratio less than 1 indicates that $\phi_1$ achieves better performance on energy efficiency than $\phi_2$) of PIER to PTR and PLPC are presented in three subfigures, where the scaling parameter $h$ is set to $1$, $10$ and $20$, respectively. Each curve is plotted based on the results of $500$ independent runs. We observe that, PIER is more efficient compared to PTR and PLPC for over $86.2\%$ and $71\%$ of simulation runs in all cases, respectively.

\subsection{Robustness to Distributions of Task Durations}

\begin{figure}[ht]
\centering
\includegraphics[height=4.5cm,width=7cm]{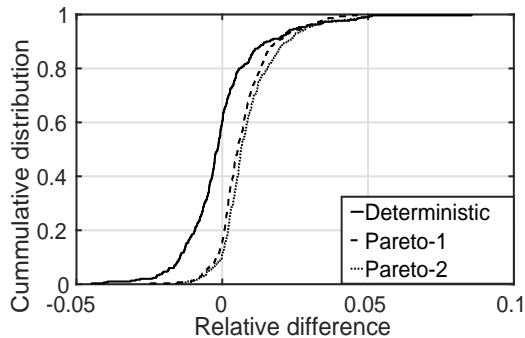}
\caption{Robustness of PIER in terms of energy efficiency under different distributions of task durations.}\label{fig:fig3}
\end{figure}

We then alter the distribution of task durations to explore the robustness of PIER.
We consider three distributions: deterministic distribution, Pareto distribution with finite variance (shape parameter = $2.002$, denoted as Pareto-1) and Pareto distribution with infinite variances (shape parameter = $1.981$, denoted as Pareto-2).
The ranges of all parameters are the same as those in Fig.~\ref{fig:fig2} with $h:=1$. In Fig.~\ref{fig:fig3}, we plot the cumulative distribution of the relative difference in power consumption between the one with the exponential distribution and the one with the specified distribution obtained over $500$ independent simulation runs under PIER. The results indicate that the mean relative differences in all tested cases are within $\pm 4\%$ for all simulation runs. That is, PIER is robust within the demonstrated cases as its performance is not very sensitive to the distribution of task durations that could be affected by nature of different offloaded tasks.

\section{Conclusions}\label{sec:conclusion}
We studied the TOSP in a heterogeneous fog computing environment including a central cloud.
We proposed a priority-based scheduling policy, PIER, which always selects the ARC group with the highest index, computed as the ratio of the mean service rate to the incremental power rate, as the prioritized offloading destination. We showed that, in a case applicable for local area wireless networks with high volumes of tasks and exponentially distributed durations, PIER is equivalent to a particular case of a policy proposed in~\cite{Fu2020} that has been proved asymptotically optimal.

In addition, we demonstrated the improvements of the proposed policy in terms of energy efficiency by extensive numerical experiments for more general cases with multiple classes of tasks and ARC groups. Compared with two benchmark policies, PTR and PLPC, that had been commonly adopted in existing work, the PIER acts more effectively in $78.6\%$ of all simulation runs, with significantly better results achieved in most  cases. We also demonstrated that, although PIER is derived under the assumption that the task durations are exponentially distributed, its robustness to the different distributions of task durations enables the possibility to apply the policy to a broader range of practical scenarios.

%
%
%
\bibliographystyle{splncs04}
\bibliography{paperbib}

\end{document}